\documentclass{article}

\usepackage{latexsym, amsfonts}
\usepackage{amssymb, amsmath}
\usepackage{amsfonts, latexsym}
\usepackage{booktabs}

\newcommand{\qed}{\hfill $\square$ \vspace{3mm}}
\newcommand{\proof}{\noindent \textbf{Proof \ }}
\newcommand{\im}{\mathrm{im\, }}
\newcommand{\rank}{\mathrm{rank\, }}

\newcommand{\fix}{\mathrm{fix\, }}

\newcommand{\opd}{\mathrm{opd\, }}
\newcommand{\ord}{\mathrm{ord\, }}

\newtheorem{theorem}{Theorem}
\newtheorem{lemma}[theorem]{Lemma}

\newtheorem{proposition}[theorem]{Proposition}

\begin{document}
	
\title{On certain subsemigroups of finite oriented and order-decreasing full transformations}

\author{Gonca Ay\i k$^{1}$,  Hayrullah Ay\i k$^{1}$, Ilinka Dimitrova$^{2}$, J\"{o}rg Koppitz$^{3,*}$\vspace{3mm}\\
		$^{1}$Department of Mathematics, \c{C}ukurova University, Adana, Turkey \\
        $^{2}$Faculty of Mathematics and Natural Science, \\
        South-West University "Neofit Rilski", Blagoevgrad, Bulgaria \\
        $^{3}$Institute of Mathematics and Informatics,\\
        Bulgarian Academy of Sciences, Sofia, Bulgaria.}

\date{}
\maketitle
\let\thefootnote\relax\footnote{$^{*}$ Corresponding author.}

\let\thefootnote\relax\footnote{\, \ E-mail adresses:  agonca@cu.edu.tr (Gonca Ay\i k), hayik@cu.edu.tr (Hayrullah Ay\i k), ilinka\_dimitrova@swu.bg (Ilinka Dimitrova), koppitz@math.bas.bg (J\"{o}rg Koppitz)}	

\let\thefootnote\relax\footnote{\, \ This study was supported by Scientific and Technological Research Council of Turkey (TUBITAK) under the Grant Number 123F463. The first and second authors thank to TUBITAK for their supports.}	

\noindent \textbf{Abstract} 
\noindent Let $\mathcal{ORD}_{n}$ be the semigroup consisting of all oriented and order-decreasing full transformations on the finite chain $X_{n}=\{ 1<\cdots<n \}$, and for $1\leq r\leq n-1$, let 
$$\mathcal{ORD}(n,r) =\{\alpha \in \mathcal{ORD}_{n}\, :\, \lvert\im(\alpha )\rvert \leq r\}.$$
In this paper, we determine the cardinality of $\mathcal{ORD}(n,r)$ and the number of nilpotent elements of $\mathcal{ORD}(n,r)$, we find a minimal generating set and the rank of $\mathcal{ORD}(n,r)$, and moreover, we characterize all maximal subsemigroups of $\mathcal{ORD}(n,r)$ for each $3\leq r\leq n-1$. 

\vspace{3mm}

\vspace{3mm}
	

\section{Introduction}

Let $\mathcal{T}_{n}$ denote the (full) transformation semigroup on the chain $X_{n}=\{1, \ldots , n\}$ under its natural order. An element $\alpha \in \mathcal{T}_{n}$ is called \emph{order-preserving} (\emph{order-reversing}) if $a<b$ implies $a\alpha\leq b\alpha$ ($b\alpha\leq a\alpha$) and \emph{order-decreasing} (\emph{extensive}) if $a\alpha\leq a$ $(a\leq a\alpha)$ for all $a,b\in X_{n}$. Furthermore, $\alpha$ is called \emph{orientation-preserving} (\emph{orientation-reversing}) if the sequence $(1\alpha ,\ldots, n\alpha)$ is cyclic (anti-cyclic), i.e. there exists no more than one element $i\in X_{n}$
such that $(i+1)\alpha <i\alpha$ ($i\alpha < (i+1)\alpha$) where $(n+1)\alpha =1\alpha$. A transformation is called \emph{oriented} (\emph{monotone}) if it is orientation-preserving or orientation-reversing (order-preserving or order-reversing). Then we denote the subsemigroups of $\mathcal{T}_{n}$ consisting of all order-preserving transformations by $\mathcal{O}_{n}$, consisting of all monotone transformations by $\mathcal{M}_{n}$, consisting of all order-decreasing (extensive) transformations by $\mathcal{D}_{n}$ ($\mathcal{E}_{n}$), consisting of all orientation-preserving transformations by $\mathcal{OP}_{n}$, and consisting of all oriented transformations by $\mathcal{OR}_{n}$. We let $\mathcal{C}_{n} =\mathcal{O}_{n} \cap \mathcal{D}_{n}$,\,  $\mathcal{OPD}_{n} =\mathcal{OP}_{n} \cap \mathcal{D}_{n}$ and $\mathcal{ORD}_{n} =\mathcal{OR}_{n} \cap \mathcal{D}_{n}$, and for $1\leq r\leq n$, we also let
\begin{eqnarray*}
	&\mathcal{U}(n,r)=\{ \alpha\in \mathcal{U}_{n} : \lvert\im(\alpha)\rvert \leq r \}&
\end{eqnarray*}
for any subsemigroup $\mathcal{U}_{n}$ of $\mathcal{T}_{n}$. If we denote the subsemigroup of $\mathcal{T}_{n}$ consisting of all oriented and extensive transformations by $\mathcal{ORE}_{n}$, then one can similarly prove as in \cite[Lemma 1.1]{U2} that $\mathcal{ORD}_{n}$ and $\mathcal{ORE}_{n}$ are isomorphic semigroups, and so it is enough to only consider the semigroup $\mathcal{ORD}_{n}$. For other terms in semigroup theory, which are not explained here, we refer to \cite{GM} and \cite{H}.

Let $S$ be a semigroup and let $A$ be a non-empty subset of $S$. Then we denote the smallest subsemigroup of $S$ containing $A$ by $\langle A\,\rangle$, and if $\langle A\, \rangle =S$, then $A$ is called a \emph{generating set} of $S$. The \emph{rank} of $S$ is defined by $\rank(S) =\min \{\lvert A\rvert :\langle A\, \rangle =S\}$, and a  generating set of $S$ with size $\rank(S)$ is called a \emph{minimal generating set} of $S$. An element $s\in S$ is called \emph{undecomposable} in $S$ if, for $u,v\in S$,\,  $s=uv$ implies $s=u$ or $s=v$, and clearly, every generating set of $S$ contains all its undecomposable elements. The set of all idempotents in a non-empty subset $T$ of $S$ is denoted by $E(T)$, i.e. $E(T)=\{ e\in T: e^{2}=e\}$. An element $a$ of a semigroup $S$ with $0$ is called \emph{nilpotent} if $a^{k}=0$ for some $k\in \mathbb{N}$, and the set of all nilpotent elements in a non-empty subset $T$ of $S$ is denoted by $N(T)$, i.e. $N(T)=\{ a\in S : a^{k}=0\, \mbox{ for some }\, k\in \mathbb{N}\}$. Furthermore, a semigroup $S$ with $0$ is called \emph{nilpotent} if $S^{m}=\{0\}$ for some $m\in\mathbb{N}$. It is known that a finite semigroup $S$ with $0$ is nilpotent if and only if every element of $S$ is nilpotent (see, for example, \cite [Proposition 8.1.2]{GM}), and it is obvious that $S\setminus S^{2}$ is a unique minimal generating set of $S$.

The problems to find (minimal) generating sets, and so ranks of a finite transformation semigroup, and to characterize maximal subsemigroups of a finite transformation semigroup have been studied for a long time in semigroup theory (see, for example, \cite{AABD, AAD, AADK, GH, U2, YK}). In more detail, it is shown in \cite[Corollary 2.4]{PH1} and \cite[Section 4]{LU2} that $\rank(\mathcal{C}_{n})= n$ and $\rank (\mathcal{C}(n,r)) =\binom{n-1}{r-1}$ for $1\leq r\leq n-1$, respectively. Moreover, it is shown in \cite[Theorem 1.4]{ZF} that $\rank(\mathcal{M}(n,r)) =\binom{n}{r}$ for $1\leq r\leq n-1$. The semigroup $\mathcal{OP}_{n}$ was introduced independently by Mc Alister in \cite{MA}, and Catarino and Higgins in \cite{CH}, and many studies have been done on $\mathcal{OP}_{n}$ and $\mathcal{OPD}_{n}$ (see, for example, \cite{FGJ, LZL, ZP}). In more detail, it is shown in \cite[Theorem 12]{D2} that $\rank(\mathcal{OPD}(n,r)) =\binom{n}{r}+r-2$ for $2\leq r\leq n-1$. Beyond all such studies, Gray developed some methods in \cite{RDG}. These methods can be used to compute the ranks of the proper two-sided ideals of $\mathcal{T}_{n}$ and its subsemigroups  generated by mappings with prescribed and images. Moreover, for $2\leq r\leq n-1$, Ay\i k and Bugay found necessary and sufficient conditions for any set of full transformations of height $r$ in $\mathcal{T}_{n}$ to be a (minimal) generating set of  $\mathcal{T}(n,r)$ in \cite{AB}. Unlike these, we examine the ranks of certain subsemigroups based on undecomposable elements. In this paper, we will find a (minimal) generating set and the rank of $\mathcal{ORD}(n,r)$ for each $3\leq r\leq n$.
 
Maximal subsemigroups were characterized for a variety of semigroups. The paper \cite{GGR} by Graham, Graham and Rhodes is basically for the study of maximal subsemigroups of finite semigroups. In \cite{DMW}, Donoven, Mitchell, and Wilson described an algorithm for calculating the maximal subsemigroups of an arbitrary finite semigroup, using the results in \cite{GGR}. East, Kumar, Mitchell, and Wilson used the framework provided by this algorithm to describe and count the maximal subsemigroups of several families of finite transformation semigroups \cite{EKMW}. Although these basic papers give a general description of maximal subsemigroups, it is important to know the maximal subsemigroups for concrete semigroups. Maximal subsemigroups of many well-known transformation semigroups have been studied in several papers (see, for example, \cite{DK, DKF, Yang, YY}). We will determine the maximal subsemigroups of $\mathcal{ORD}(n,r)$ for each $3\leq r\leq n$.

\section{Cardinalities}

The \emph{fix} and \emph{kernel} sets of $\alpha\in \mathcal{T}_{n}$ are defined by
\begin{eqnarray*}
	\fix(\alpha) &=& \{ a\in X_{n} : a\alpha=a\} \,  \mbox{ and }\\
	\ker(\alpha) &=& \{ (a,b)\in X_{n}\times X_{n} : a\alpha= b\alpha\},
\end{eqnarray*} 
respectively. It is clear that $\alpha\in \mathcal{T}_{n}$ is an idempotent if and only if $\fix(\alpha) =\im(\alpha)$, and that $\alpha\in \mathcal{D}_{n}$ is an idempotent if and only if $\min(a\alpha^{-1}) =a$ for each $a\in \im(\alpha)$. For a non-empty subset $Y$ of $X_{n}$, the restriction map of $\alpha\in \mathcal{T}_{n}$ to $Y$ is denoted by $\alpha_{\mid_{Y}}$, and moreover, $\alpha_{\mid_{Y}}$ is call the \emph{partial identity} on $Y$ if $a\alpha_{\mid_{Y}} =a$ for all $a\in Y$. However, the identity transformation on $X_{n}$ is denoted by $1_{n}$. 

For any two non-empty subsets $A$ and $B$ of $X_{n}$, we write $A<B$ if $a<b$ for all $a\in A$ and $b\in B$. A non-empty subset $C$ of $X_{n}$ is said to be a \emph{convex} subset of $X_{n}$ if $a\leq c\leq b$ implies $c\in C$ for all $a,b \in C$ and $c\in X_{n}$. Let $P=\{ A_{1}, \ldots ,A_{r} \}$ ($1\leq r\leq n$) be a partition of $X_{n}$ into $r$ subsets, that is, a family of disjoint non-empty $r$ subsets of $X_{n}$ whose union is $X_{n}$. If $A_{i}< A_{i+1}$ for all $1\leq i\leq r-1$, then $P$ is called  an \emph{ordered} partition of $X_{n}$, and we write $P=(A_{1}< \cdots <A_{r})$. If each $A_{i}$ is convex for $1\leq i\leq r$, then $P$ is called a \emph{convex} partition of $X_{n}$. Notice that for each $\alpha\in \mathcal{T}_{n}$,\, $\ker(\alpha)$ is an equivalence relation on $X_{n}$, and the set of equivalence classes, namely $\{ a\alpha^{-1}: a\in \im(\alpha)\}$, is a partition of $X_{n}$. In addition, a subset $Y$ of $X_{n}$ is called a \emph{transversal} of $\ker(\alpha)$ if $\lvert Y\cap a\alpha^{-1}\rvert =1$ for all $a\in \im(\alpha)$. Notice also that for $\alpha\in \mathcal{O}_{n}$, if $\im(\alpha) =\{ a_{1}< \cdots< a_{r}\}$, then $(a_{1}\alpha^{-1} <\cdots< a_{r}\alpha^{-1})$ is a convex and ordered partition of $X_{n}$, and for $\alpha\in \mathcal{M}_{n} \setminus \mathcal{O}_{n}$, if $\im(\alpha) =\{ a_{1}< \cdots< a_{r}\}$, then $(a_{r}\alpha^{-1} <\cdots< a_{1}\alpha^{-1})$ is a convex and ordered partition of $X_{n}$. Moreover, it is clear that  
$\pi =\left( \begin{matrix}
	1&2&\cdots &n\\
	1&1&\cdots &1
\end{matrix}\right)$ is the zero element of $\mathcal{D}_{n}$, and so of $\mathcal{ORD}_{n}$. For any $\alpha \in \mathcal{D}_{n}$, since $\alpha^{k} =\pi$ for some $k\in \mathbb{N}$ is equivalent to $\fix(\alpha)= \{1\}$, it follows that $N(\mathcal{D}_{n})=\{ \alpha\in \mathcal{D}_{n} : \fix (\alpha )=\{ 1\} \}$, and so $N(\mathcal{ORD}_{n})=\{ \alpha\in \mathcal{ORD}_{n} : \fix (\alpha )=\{ 1\} \}$. 

Let $\mathcal{R}_{n}$ denote the subset of $\mathcal{T}_{n}$ consisting of all orientation-reversing transformations, and let $\mathcal{RD}_{n} =\mathcal{R}_{n} \cap \mathcal{D}_{n}$. It is shown in \cite[Lemma 1.1]{CH} that for $\alpha\in \mathcal{T}_{n}$, the sequence $(1\alpha , \ldots ,n\alpha)$ is both cyclic and anti-cyclic if and only if $\lvert \im(\alpha) \rvert \leq 2$. Thus, we have 
$$\mathcal{OP}_{n} \cap \mathcal{R}_{n} =\{ \alpha \in \mathcal{OR}_{n} : \lvert \im(\alpha) \rvert \leq 2 \}.$$ 
It is shown in \cite[Proposition 2.3]{CH} that a non-constant $\alpha$ in $\mathcal{OR}_{n}$ is order-preserving (order-reversing) if and only if $1\alpha < n\alpha$ ($n\alpha < 1\alpha$). Since $\mathcal{ORD}_{n}= \mathcal{OPD}_{n} =\mathcal{C}_{n}$ for $n=1,2$, and $\mathcal{ORD}_{3} = \mathcal{OPD}_{3}$, we suppose $n\geq 4$ throughout this paper unless otherwise stated. Since $\mathcal{ORD}(n,1)= \mathcal{OPD}(n,1)= \mathcal{C}(n,1)$,  $\mathcal{ORD}(n,2)= \mathcal{OPD}(n,2)$, $\mathcal{ORD}(n,n-1)= \mathcal{ORD}_{n} \setminus \{ 1_{n}\}$ and $\mathcal{ORD}(n,n)= \mathcal{ORD}_{n}$, we also suppose $3\leq r\leq n-1$ throughout this paper unless otherwise stated.  Next we let 
$$\mathcal{RD}_{n}^{*} =\mathcal{RD}_{n} \setminus \mathcal{OPD}_{n} =\mathcal{ORD}_{n} \setminus \mathcal{OPD}_{n},$$
and for $1\leq i\leq j\leq n$, we put 
$$[i,j]=\{ i,\ldots, j\}.$$
For each $\alpha\in \mathcal{RD}_{n}^{*}$, it is clear that $\rvert \im(\alpha) \lvert \geq 3$, and moreover, if $\im(\alpha) =\{1< a_{2}< \cdots< a_{r}\}$, then it is also clear that $(a_{r}\alpha^{-1}, \ldots, a_{2}\alpha^{-1})$ is a convex and ordered partition of $X_{n}\setminus (1\alpha^{-1})$. Thus, we conclude that   
\begin{eqnarray}\label{e1}
r\leq a_{r} \leq \min(a_{r}\alpha^{-1})= \min(X_{n} \setminus (1\alpha^{-1}))&\mbox{and}&
r\leq \left \lceil \frac{n+1}{2} \right \rceil,
\end{eqnarray}
where $\left \lceil \frac{n+1}{2} \right \rceil$ denotes the smallest integer not strictly less than $\frac{n+1}{2}$. 

In \cite[Section 2]{D2}, the \emph{order-preserving degree} of $\alpha\in \mathcal{OPD}_{n}$ is defined by $\opd (\alpha) =\max \{ m : \alpha_{\mid_{X_{m}}}\in  \mathcal{O}_{m}\}$. Then noticed that for $\alpha\in \mathcal{OPD}_{n} \setminus \mathcal{C}_{n}$, if $\opd(\alpha) =m$, then $2\leq m\leq n-1$ and $\alpha$ has the following tabular form:
\begin{eqnarray*}
	\alpha=\left(
	\begin{array}{cccc|cccc}
		1 & \cdots & m-1 & m &  m+1 & \cdots & n\\
		1 & \cdots & (m-1)\alpha & m\alpha &  1 & \cdots & 1 \\
	\end{array}\right) 
\end{eqnarray*}
with the property that $1\leq 2\alpha\leq \cdots\leq (m-1)\alpha\leq m\alpha\leq m$ and $m\alpha \geq 2$.

For $\alpha\in \mathcal{RD}_{n}^{*}$, we now define the \emph{order-reserving degree} of $\alpha$ by 
$$\ord (\alpha) =\min(X_{n} \setminus (1\alpha^{-1})).$$ 
Then, for $\alpha\in \mathcal{RD}_{n}^{*}$, if $\ord(\alpha) =m$, then we notice that $3\leq m \leq n-1$ and $\alpha$ has the following tabular form:
\begin{eqnarray}\label{e2}
\alpha=\left(\begin{array}{ccc|ccccc}
	1 & \cdots & m-1 &    m    &     m+1     & \cdots & n\\
	1 & \cdots &  1  & m\alpha & (m+1)\alpha & \cdots & n\alpha \\
\end{array}\right)
\end{eqnarray}
with the property that $1\leq n\alpha\leq \cdots\leq (m+1)\alpha\leq m\alpha\leq m$ and $m\alpha \geq 3$. In addition, if $\lvert \im(\alpha) \rvert=r$, then we also notice that 
\begin{eqnarray*}
	r\leq \min\left\{ m,\left \lceil\frac{n+1}{2} \right \rceil \right\}\,\, \mbox{ and }\,\, \fix(\alpha) \subseteq \{1,m\}.
\end{eqnarray*}
For each $r\leq \left \lceil\frac{n+1}{2} \right \rceil$, we let
$$J_{n,r}=\{ \alpha\in \mathcal{RD}_{n}^{*} :\lvert \im(\alpha) \rvert=r\}\, \mbox{ and we have }\, N(J_{n,r})=\{ \alpha\in J_{n,r} : \fix(\alpha) =\{1\}\}.$$
Before giving our first result, we want to remind the Rothe-Hagen identity, which is a generalization of Vandermonde's identity:
\begin{eqnarray*}
	\sum\limits_{k=0}^{m} \frac{a}{a+bk}\binom{a+bk}{k} \binom{c-bk}{m-k}= \binom{a+c}{m}
\end{eqnarray*}
for all $a,b,c\in \mathbb{C}$ (see, for example \cite[Section 5.4]{GKP}). If we take $a=r$,\, $b=1$,\ $c=n-r+1$ and $m=n-2r+2$, then we have
\begin{eqnarray*}
	\sum\limits_{k=0}^{n-2r+2} \binom{r-1+k}{k} \binom{n-r+1-k}{n-2r+2-k}= \binom{n+1}{n-2r+2}=\binom{n+1}{2r-1}.
\end{eqnarray*}

\begin{proposition}\label{p1}
	Let $3\leq r\leq \left \lceil\frac{n+1}{2} \right \rceil$. Then we have $\lvert J_{n,r}\rvert =\binom{n+1}{2r-1}$, and moreover, $\lvert N(J_{n,r}) \rvert= \lvert J_{n-1,r}\rvert$ for $n\geq 5$.
\end{proposition}

\proof Let $G(n,r,m)=\{ \alpha\in J_{n,r} : \ord(\alpha)=m \}$, and let $\alpha\in G(n,r,m)$. Then $x\alpha =1$ for all $x\in X_{m-1}$ and $\alpha$ is an order-reversing mapping from $[m,n]$ to $X_{m}$ with the property that $(\im(\alpha) \setminus\{1\}) \subseteq \im(\alpha_{\mid_{[m,n]}})$ and $3\leq m\alpha =\max(\im(\alpha)) \leq m$. More precisely, if $\im(\alpha) =\{1<a_{2} <\cdots< a_{r}\}$, then $r\leq m\leq n-r+2$ and $\alpha$ has one of the following tabular forms:
$$\alpha=  \left( \begin{array}{cccccccc}
	X_{m-1} & A_{r} &\cdots & A_{2} \\
	   1    & a_{r} &\cdots & a_{2}
\end{array}\right)\, \mbox{ or }\, 
\alpha=  \left( \begin{array}{cccccccc}
	X_{m-1} & B_{r} &\cdots & B_{2} & B_{1} \\
	   1    & a_{r} &\cdots & a_{2} &   1
\end{array}\right)$$
where $(A_{r}< \cdots< A_{2})$ and $(B_{r}< \cdots< B_{2}< B_{1})$ are convex and ordered partitions of $[m,n]$. Therefore, since there exist $\binom{m-1}{r-1}$ subsets of $X_{m}\setminus \{1\}$ with size $r-1$, and since there exist $\binom{n-m}{r-2}$ convex and ordered partitions of $[m,n]$ into $r-1$ subsets and $\binom{n-m}{r-1}$ convex and ordered partitions of $[m,n]$ into $r$ subsets, it follows that
$$\lvert G(n,r,m) \rvert =\binom{m-1}{r-1}\left(\binom{n-m}{r-2} + \binom{n-m}{r-1} \right) =\binom{m-1}{r-1}\binom{n-m+1}{r-1}.$$
Finally, since $r\leq m\leq n-r+2$, we have 
\begin{eqnarray*}
	\lvert J_{n,r}\rvert &=&\sum\limits_{m=r}^{n-r+2} \binom{m-1}{r-1} \binom{n-m+1}{r-1} =\sum\limits_{k=0}^{n-2r+2} \binom{r-1+k}{r-1} \binom{n-r+1-k}{r-1}\\
	&=&\sum\limits_{k=0}^{n-2r+2} \binom{r-1+k}{r-1} \binom{n-r+1-k}{n-2r+2-k}= \binom{n+1}{2r-1}.
\end{eqnarray*}

First of all, it is clear that $N(\mathcal{RD}_{4}^{*}) =\emptyset$. Suppose that $n\geq 5$ and for each $\alpha\in \mathcal{RD}_{n-1}^{*}$, define $\hat{\alpha} :X_{n} \rightarrow X_{n}$ by $1\hat{\alpha}=1$ and $i\hat{\alpha} =(i-1)\alpha$ for all $2\leq i\leq n$. Then it is clear that $\hat{\alpha} \in N(\mathcal{RD}_{n}^{*})$ and $\lvert \im(\alpha) \rvert= \lvert \im(\hat{\alpha}) \rvert$. Moreover, consider the mapping $\psi :\mathcal{RD}_{n-1}^{*} \rightarrow N(\mathcal{RD}_{n}^{*})$ defined by $\alpha \psi= \hat{\alpha}$ for all $\alpha\in \mathcal{RD}_{n-1}^{*}$. It is also clear that $\psi$ is a bijection. Therefore, we have $\lvert N(J_{n,r}) \rvert= \lvert J_{n-1,r} \rvert$, as required. \qed

Next, we determine the cardinalities of $\mathcal{RD}_{n}^{*}$ and $N(\mathcal{RD}_{n}^{*})$ in the following lemma.

\begin{lemma}\label{l2} We have
	\begin{itemize}
		\item[$(i)$] $\lvert \mathcal{RD}_{n}^{*} \rvert= 2^{n}-\frac{n^{3}+5n+6}{6}$ and
		\item[$(ii)$] $\lvert N(\mathcal{RD}_{n}^{*}) \rvert= \lvert \mathcal{RD}_{n-1}^{*} \rvert$ for all $n\geq 5$.
	\end{itemize}
\end{lemma}

\proof $(i)$ Let $\alpha\in \mathcal{RD}_{n}^{*}$. If $\ord(\alpha)=m$, then, as observed above, $3\leq m\leq n-1$,\, $\lvert \im(\alpha) \rvert \geq 3$, and $\alpha$ has the tabular form as defined in (\ref{e2}).

For each $1\leq i\leq m$, let $c_{i}=\lvert \{ a\in [m,n] : a\alpha =i\}\rvert$. Then it is clear that $(c_{1}, c_{2}, \ldots, c_{m})$ is a solution of the equation:
$$x_{1}+x_{2}+\cdots+x_{m}=n-m+1\, \mbox{ and }\, x_{1},x_{2},\ldots,x_{m} \geq 0.$$
Moreover, for each solution of this equation, one can easily define a mapping in $\mathcal{RD}_{n}^{*}$. The number of ordered integers solutions of the above equation is $\binom{n}{m-1}$. Clearly, some solutions give mappings in $\mathcal{OPD}_{n}$, i.e. mappings in $\mathcal{RD}_{n}$ whose image sets contain at most $2$ elements. Since there exist $m$ mappings whose restriction to $[m,n]$ are constant, and since for all $2\leq k\leq n$ there exist $(m-1)(n-m)$ mappings whose restricted images to $[m,n]$ are $\{1,k\}$, it follows that
$$\lvert\{ \alpha\in \mathcal{RD}_{n}^{*} :\ord(\alpha) =m\}\rvert =\binom{n}{m-1} -m-(m-1)(n-m).$$ 
Therefore, we have
\begin{eqnarray*}
\lvert \mathcal{RD}_{n}^{*} \rvert= \sum\limits_{m=3}^{n-1} \left(\binom{n}{m-1} +m^{2}-m(n+2)+n \right)= 2^{n}-\frac{n^{3} +5n+6}{6}.
\end{eqnarray*}

$(ii)$ By the bijection defined in the proof of Proposition \ref{p1}, we also have $\lvert N(\mathcal{RD}_{n}^{*}) \rvert= \lvert \mathcal{RD}_{n-1}^{*} \rvert$ for all $n\geq 5$.  \qed

Let $n\geq 1$ and $1\leq r\leq n$. Then $n^{\text{th}}$ \emph{Catalan number} $C_{n}$  and \emph{Narayana number} $N(n,r)$ are defined by $C_{n}= \frac{1}{n+1} \binom{2n}{n}$ and $N(n,r)= \frac{1}{n} \binom{n}{r} \binom{n}{r-1}$, respectively. For each $1\leq r\leq n$, it is shown in \cite[Lemma 1]{D2} that 
$$\lvert\mathcal{OPD}(n,r)\rvert =-n+1+\sum\limits_{m=1}^{r} C_{m}+ \sum\limits_{m=r+1}^{n} \sum\limits_{k=1}^{r} N(m,k),$$ 
and that $\lvert N(\mathcal{OPD}(n,r))\rvert = \lvert \mathcal{OPD}(n-1,r)\rvert$. Therefore, we are able to state the main result of this section.

\begin{theorem}  We have
	\begin{itemize}
		\item[$(i)$] $\lvert \mathcal{ORD}(n,r)\rvert = \lvert \mathcal{OPD}(n,r)\rvert+ \sum\limits_{k=3}^{\hat{r}} \binom{n+1}{2k-1}$, where $\hat{r} = \left\{\begin{array}{cl}
			r & \textrm{if } r\leq \left \lceil\frac{n+1}{2} \right \rceil ,\vspace{1mm}\\
			\left \lceil\frac{n+1}{2} \right \rceil & \textrm{if } r > \left \lceil\frac{n+1}{2} \right \rceil
		\end{array}\right.$ and
		\item[$(ii)$] $\lvert N(\mathcal{ORD}(n,r))\rvert = \lvert \mathcal{ORD}(n-1,r)\rvert $ for $n\geq 5$.
	\end{itemize}
Furthermore, $\lvert \mathcal{ORD}_{n}\rvert = 2^{n}-\frac{n^{3}+11n}{6}+\sum\limits_{m=1}^{n} C_{m}$.  \hfill $\square$
\end{theorem}

\section{Generating set and rank of $\mathcal{ORD}(n,r)$}

After the experience we gained in the previous section, it is now clear that $E(\mathcal{ORD}_{n}) =E(\mathcal{OPD}_{n})$. For each $1\leq r\leq n$ and any subset $\mathcal{S}$ of $\mathcal{T}_{n}$, we let 
$$E_{r}(\mathcal{S})=\{ \alpha \in E(\mathcal{S}) : \lvert\im(\alpha) \rvert =r \}.$$
For each subset $Y$ of $X_{n}$ which contains $1$, it is also clear that there exists a unique idempotent in $\mathcal{C}_{n}$ whose image is $Y$, namely if $Y=\{1< a_{2}< \cdots< a_{r}\}$, then 
$\zeta_{Y}=\left( \begin{smallmatrix}
	[1,a_{2}-1] & [a_{2},a_{3}-1] &\cdots & [a_{r-1},a_{r}-1] & [a_{r},n]\\
	     1      &      a_{2}      &\cdots &      a_{r-1}      &   a_{r}
\end{smallmatrix}\right)$
(see, also \cite[Lemma 3.12]{LU2}).  It is shown in \cite{LU2} that $\fix(\alpha \beta) =\fix(\alpha)\cap \fix(\beta)$ for all $\alpha , \beta \in \mathcal{D}_{n}$, and that $E_{r}(\mathcal{C}_{n})$ is a unique minimal generating set of $\mathcal{C} (n,r)$, and $\rank (\mathcal{C}(n,r))= \binom{n-1}{r-1}$. For each $1\leq r\leq n-1$ and $r+1\leq m\leq n$, we let 
$$B_{m}=\{ \alpha\in E_{r}(\mathcal{OPD}_{n}) : \opd(\alpha)=m \}.$$ 
Then we immediately have $\lvert B_{m} \rvert =\binom{m-1}{r-1}$. Also let 
$$C=\left(\bigcup\limits_{m=r+1}^{n}B_{m}\right)\cup \{\rho_{2},\ldots, \rho_{r}\},$$
where $\rho_{m}=\left( \begin{matrix}
	1 & 2 &\cdots & m & m+1 &\cdots & n\\
	1 & 2 &\cdots & m &  1  &\cdots & 1
\end{matrix}\right)$ for each $2\leq m\leq r$. Then we have the following result from \cite{D2}.

\begin{theorem}\label{t3} \cite[Theorem 12]{D2}
	For $2\leq r\leq n-1$,\, $C$ is the unique minimal generating set of $\mathcal{OPD}(n,r)$, and so $\rank (\mathcal{OPD}(n,r))=\binom{n}{r}+r-2$. \hfill $\square$
\end{theorem}

Next we state and prove the following result:

\begin{proposition}\label{p4}
	For $2\leq r\leq n-1$, each element of $C$ is undecomposable in  $\mathcal{ORD}(n,r)$.
\end{proposition}

\proof Let $2\leq r\leq n-1$ and let $\zeta$ be any element in $C$ with $\im(\zeta) =\{ 1=a_{1}< a_{2}< \cdots <a_{r}\}$. Assume that $\zeta =\alpha \beta$ for some $\alpha, \beta\in \mathcal{ORD}(n,r)$. Since $r=\lvert \im(\zeta) \rvert \leq \lvert \im(\alpha) \rvert$, we first have $\lvert \im(\alpha) \rvert= r$, and moreover, since $\ker(\alpha) \subseteq \ker(\zeta)$, we have $\ker(\zeta)= \ker(\alpha)$. Thus, since $\fix(\zeta) \subseteq \fix(\alpha)$, we have $a_{i}\alpha =a_{i}$, and so $(a_{i}\zeta^{-1})\alpha =\{a_{i}\}$ for each $1\leq i\leq r$. Thus, we have $\zeta=\alpha$. 

Next let $2\leq m\leq r-1$, and assume that $\rho_{m} =\alpha \beta$ for some $\alpha, \beta\in \mathcal{ORD}(n,r)$. Because of order-decreasing, we notice that $\fix(\rho_{m}) \subseteq \fix(\alpha) ,\, \fix(\beta)$, and hence we have $i\alpha =i$ and $i\beta =i$ for all $1\leq i\leq m$, and moreover, we have $(m+1)\alpha \leq m+1$. If $(m+1)\alpha =m+1$, then it is clear that $\rho_{m} =\beta$. If $(m+1)\alpha \leq m$, then we have $(m+1)\alpha= ((m+1)\alpha) \beta =(m+1)\rho_{m}  =1$. Thus, we have $\rho_{m} =\alpha$. \qed

For each $3\leq r\leq \left \lceil\frac{n+1}{2} \right \rceil$, it is easy to verify that   
$$\mathcal{RD}^{*}(n,r) =\mathcal{ORD}(n,r) \setminus \mathcal{OPD}(n,r).$$
Next, for each $3\leq r\leq \left \lceil\frac{n+1}{2} \right \rceil$ and $3\leq m\leq n-1$, we define $\lambda_{m,r}$ in $\mathcal{RD}^{*}(n,r)$ as follows: If $3\leq m\leq n-r+1$, then 
\begin{eqnarray}\label{e4}
\lambda_{m,r}=\left(\begin{array}{cccccccccc}
		[1,m-1] & m & m+1 &\cdots & 2m-p-2 & [2m-p-1,n]\\
		1    & m & m-1 &\cdots &   p+2  &    1
\end{array}\right),
\end{eqnarray}
where $p=\max\{0,m-r\}$. If $n-r+2\leq m\leq n-1$, then 
\begin{eqnarray}\label{e5}
	\lambda_{m,r}=\left(\begin{array}{cccccccccc}
		[1,m-1] & m & m+1 &\cdots &   n  \\
		1    & m & m-1 &\cdots & 2m-n 
	\end{array}\right).
\end{eqnarray}
Notice that $r\leq \left \lceil\frac{n+1}{2} \right \rceil$ implies $r \leq \frac{n}{2} + 1$ and thus $2r-2 \leq n$. Hence, $r=2r-r-2+2 \leq n-r+2$. So $n-r+2 \leq m$ implies $r \leq m$ and $2\leq m+r-n \leq 2m-n$. Further, for each $3\leq r\leq n-1$, we put  
$$G_{n,\hat{r}}=\{\lambda_{m,\hat{r}} : 3\leq m\leq n-1\}\,
\mbox{ where }\, \hat{r} = \left\{\begin{array}{cl}
	r & \textrm{if } r\leq \left \lceil\frac{n+1}{2} \right \rceil ,\vspace{1mm}\\
	\left \lceil\frac{n+1}{2} \right \rceil & \textrm{if } r > \left \lceil\frac{n+1}{2} \right \rceil .
\end{array}\right.$$
Then, we state and prove the following result.

\begin{lemma}\label{l5}
	$\mathcal{ORD}(n,r) =\langle \mathcal{OPD}(n,r) \cup G_{n,\hat{r}} \rangle$.
\end{lemma}

\proof It is enough to show that $\mathcal{RD}^{*}(n,r) \subseteq \langle \mathcal{OPD}(n,r) \cup G_{n,\hat{r}}\rangle$. Let $\alpha\in \mathcal{RD}^{*}(n,r)$ with $\im(\alpha) =\{1< a_{2}< \cdots< a_{k}\}$ ($3\leq k\leq r$), and let $\ord(\alpha) =m$. Moreover, we denote $a_{i}\alpha^{-1}$ by $A_{i}$ for $2\leq i\leq k$, and we let $A_{1}= 1\alpha^{-1} \cap [m,n]$ (it is possible that $A_{1}=\emptyset$). Then $\alpha$ can be written in the following tabular form:
$$\alpha =\left(\begin{array}{cccccc}
	X_{m-1} & A_{k} & A_{k-1} &\cdots & A_{2} & A_{1}\\
	1    & a_{k} & a_{k-1} &\cdots & a_{2} &  1
\end{array}\right).$$
Since $2\leq a_{2}< \cdots <a_{k}\leq m$, $A_{k}< \cdots< A_{2}$ and $A_{k}< \cdots< A_{2}< A_{1}$ if $A_{1}\neq \emptyset$, it is clear that $a_{i}\leq m-k+i$ and $m+k-i\leq \min(A_{i})$ for each $2\leq i\leq k$, and that
$m+k-1\leq \min(A_{1})$ if $A_{1}\neq \emptyset$. Thus, if we define
$$\beta_{1} =\left\{ \begin{array}{cl}
	\left(\begin{array}{ccccc}
		X_{m-1} & A_{k} & A_{k-1} &\cdots & A_{2} \\
	       1    &   m   &   m+1   &\cdots & m+k-2 
	\end{array}\right) & \mbox{if } A_{1}= \emptyset \vspace{2mm} \\
	\left(\begin{array}{cccccc}
		X_{m-1} & A_{k} & A_{k-1} &\cdots & A_{2} & A_{1}\\
		   1    &   m   &   m+1   &\cdots & m+k-2 &   1
	\end{array}\right) & \mbox{if } A_{1}\neq \emptyset 
\end{array}\right. $$ 
and
$$\beta_{2} =\left(\begin{array}{ccccccccc}
	[1,m-k+1] & m-k+2 &\cdots &   m-1   & [m,n]\\
	1     & a_{2} &\cdots & a_{k-1} & a_{k}
\end{array}\right),$$
then it is easy to see that $\beta_{1}, \beta_{2} \in \mathcal{OPD}(n,r)$. If $m+k-2=n$, then $n-r+2=m+k-r\leq m\leq n-1$ and $2m-n=m-k+2$. Thus, we have $\alpha =\beta_{1} \lambda_{m,\hat{r}} \beta_{2}$, where $\lambda_{m,\hat{r}} \in G_{n,\hat{r}}$ is defined in (\ref{e5}). 

Furthermore, if $m+k-2\leq n-1$, then we define 
$$\gamma_{m,k}=\left(\begin{array}{cccccccccc}
	[1,m-1] & m & m+1 &\cdots & m+k-2 & [m+k-1,n]\\
	1    & m & m-1 &\cdots & m-k+2 &     1 
\end{array}\right)$$
and 
$$\delta_{m,k}=\left(\begin{array}{cccccccccc}
	[1,m-k+1] & m-k+2 &\cdots & m-1 & [m,n]\\
	1     & m-k+2 &\cdots & m-1 &   m
\end{array}\right).$$ 
It is also easy to see that $\gamma_{m,k} \in \mathcal{RD}^{*}(n,r)$ and $\delta_{m,k} \in \mathcal{C}(n,r)$. Since $m<n-k+2 \leq n-r+2$, we have $\gamma_{m,k} =\lambda_{m,\hat{r}} \delta_{m,k}$, where $\lambda_{m,\hat{r}} \in G_{n,\hat{r}}$ is defined in (\ref{e4}). Thus, we have $\alpha =\beta_{1} \gamma_{m,k} \beta_{2} =\beta_{1} \lambda_{m,\hat{r}} \delta_{m,k} \beta_{2}$ which completes the proof. \qed 

Let $\alpha$ be any element in $\mathcal{RD}^{*}(n,r)$ with $\ord(\alpha)=m$. Then notice that $\fix(\alpha) \subseteq \{1,m\}$ for some $3\leq m\leq n-1$. Next, we have the following proposition that we need for the rank of $\mathcal{ORD}(n,r)$.

\begin{proposition}\label{p6}
	Every generating set of $\mathcal{ORD}(n,r)$ must contain at least one element $\alpha \in \mathcal{RD}^{*}(n,r)$ such that $\fix(\alpha )=\{1,m\}$ for each $3\leq m\leq n-1$.
\end{proposition}

\proof Let $A$ be a generating set of $\mathcal{ORD}(n,r)$, and consider $\lambda_{m, \hat{r}}$ in $G_{n,\hat{r}}$ for each $3\leq m\leq n-1$. Then there exist $\alpha_{1}, \ldots, \alpha_{t} \in A$ such that $\lambda_{m,\hat{r}} =\alpha_{1} \cdots  \alpha_{t}$. Since $\lambda_{m,\hat{r}}$, $\alpha_{1}$, \ldots,  $\alpha_{t}$ are all order-decreasing, we have $\fix(\lambda_{m,\hat{r}}) \subseteq \fix(\alpha_{i})$, i.e. $\{1,m\} \subseteq \fix(\alpha_{i})$ for each $1\leq i\leq t$. Since the product of orientation-preserving mappings is orientation-preserving, there exists at least one $1\leq i\leq t$ such that $\alpha_{i}$ is orientation-reversing, and so $\fix(\alpha_{i})=\{1,m\}$, as required. \qed

Notice that $1_{n}$ is undecomposable in $\mathcal{ORD}_{n}$ since $\mathcal{ORD}_{n} =\mathcal{ORD}(n,n-1) \cup \{1_{n}\}$. It follows from Lemma \ref{l5} and Theorem \ref{t3} that $C\cup G_{n,\hat{r}}$ is a generating set of $\mathcal{ORD}(n,r)$. Furthermore, by Propositions \ref{p4} and \ref{p6}, we observe that $C\cup G_{n,\hat{r}}$ is minimal, and so we are able to state one of the main theorems of this paper.

\begin{theorem}\label{t7} 
	The set $C\cup G_{n,\hat{r}}$ is a minimal generating set of $\mathcal{ORD}(n,r)$, and so $\rank (\mathcal{ORD}(n,r))=\binom{n}{r}+n+r-5$. Therefore, we have $\rank (\mathcal{ORD}_{n})=3n-5$.  \hfill $\square$
\end{theorem}

Notice that $\mathcal{OPD}(n,r)$ is a semiband, that is, it is generated by a subset of $E(\mathcal{OPD}(n,r))$, for all $1\leq r\leq n$. However, $\mathcal{ORD}(n,r)$ is not a semiband since $E(\mathcal{RD}^{*}(n,r))=\emptyset$ for all $3\leq r\leq n$.

\section{Maximal Subsemigroups}

We first determine which $\lambda_{m,\hat{r}}$'s are undecomposable with the notations introduced in the previous section. 

Let $3\leq m\leq n-\hat{r}+1$. Then we define 
\begin{eqnarray*}
	\alpha_{m,\hat{r}}=\left(\begin{array}{ccccccccc}
		[1,m-1] & m & m+1 &\cdots & 2m-p-2 & [2m-p-1,n]\\
		   1    & m & m+1 &\cdots & 2m-p-2 &    1
	\end{array}\right)
\end{eqnarray*} 
and
\begin{eqnarray*}
	\beta_{m,\hat{r}}=\left(\begin{array}{cccccccccc}
		[1,m-1] & m & m+1 &\cdots & 2m-p-3 & [2m-p-2,n]\\
		   1    & m & m-1 &\cdots &   p+3  &    p+2
	\end{array}\right), 
\end{eqnarray*}
where $p=\max\{0,m-\hat{r}\}$. First, it is clear that $\alpha_{m,\hat{r}}, \beta_{m,\hat{r}} \in \mathcal{ORD}(n,r)$ and that $\alpha_{m,\hat{r}} \beta_{m,\hat{r}} = \lambda_{m,\hat{r}}$. Thus, for each $3\leq m\leq n-r+1$, we conclude that $\lambda_{m,\hat{r}}$, as defined in (\ref{e4}), is not undecomposable in $\mathcal{ORD}(n,r)$. However, if $n-\hat{r}+2 \leq m \leq n-1$, we have the following result for $\lambda_{m,\hat{r}}$'s as defined in (\ref{e5}).

\begin{proposition}\label{p8}
	Let $n-\hat{r}+2\leq m\leq n-1$. Then each $\lambda_{m,\hat{r}}$ is undecomposable in $\mathcal{ORD}(n,r)$.
\end{proposition}

\proof Let $\alpha, \beta\in \mathcal{ORD}(n,r)$ such that $\alpha\beta =\lambda_{m,\hat{r}}$. Note that, $1,m \in \fix(\lambda_{m,\hat{r}})$ implies $1,m\in \fix(\alpha) \cap \fix(\beta)$ since both $\alpha$ and $\beta$ are order-decreasing. Then since $\lambda_{m,\hat{r} \mid_{[m,n]}}$ is injective and monotone (order-reversing), it follows that $\alpha_{ \mid_{[m,n]}}$ is also injective and monotone. If $\alpha_{\mid_{[m,n]}}$ is injective and order-preserving, then $m\alpha=m$ implies that $\alpha_{\mid_{[m,n]}}$ is a partial identity, and so we obtain $\beta_{\mid_{[m,n]}} =\lambda_{m,\hat{r}\mid_{[m,n]}}$. Since $1\beta=1$, $m\beta >(m+1)\beta$, and $\beta$ is oriented, we can conclude that $[1,m-1]\beta \cap [2,m-1]=\emptyset$. Since $\beta$ is order-decreasing, we get $[1,m-1]\beta \subseteq [1,m-1]$, i.e. $[1,m-1]\beta =\{1\}$. This shows $\beta=\lambda_{m,\hat{r}}$. 

Suppose now that $\alpha_{\mid_{[m,n]}}$ is order-reversing. Since $\beta$ is order-decreasing and $\im(\alpha \beta_{\mid_{[m,n]}}) =\im(\lambda_{m, \hat{r} \mid_{[m,n]}}) =[2m-n,m]$, we obtain $\im(\alpha_{\mid_{[m,n]}}) =[2m-n,m]$, i.e. $\alpha_{\mid_{[m,n]}}$ is an order-reversing bijection from $[m,n]$ onto $[2m-n,m]$. This shows  $\alpha_{\mid_{[m,n]}} =\lambda_{m,\hat{r}\mid_{[m,n]}}$ and by the same arguments as above, we obtain $\alpha=\lambda_{m,\hat{r}}$, as required. \qed

Now, we let $3\leq m\leq n-1$ and prove some results for the set 
\begin{eqnarray*}
	L_{m,\hat{r}} &=& \{\alpha \in \mathcal{ORD}_n :\, \im(\alpha) = \im(\lambda_{m,\hat{r}}) \textrm{ and } \min(x\alpha^{-1}) = \min(x\lambda_{m,\hat{r}}^{-1})\\
	&&\qquad \qquad \qquad \textrm{ for all } x\in \im(\lambda_{m,\hat{r}})\}
\end{eqnarray*} 
that we need for the characterization of the maximal subsemigroups of $\mathcal{ORD}(n,r)$.

\begin{lemma}\label{l9}
   Let $3\leq m\leq n-\hat{r}+1$ and $\alpha \in L_{m,\hat{r}}$. Then $L_{m,\hat{r}} \subseteq \left\langle \mathcal{OPD}(n,r),\alpha \right\rangle$.
\end{lemma}

\proof 
Let $\beta \in L_{m,\hat{r}}$. Then $\im(\beta)= \im(\alpha)= \im(\lambda_{m, \hat{r}}) =\{1\}\cup [p+2,m]$, where $p=\max\{0, m-\hat{r}\}$, and $\min(x\beta^{-1})= \min(x\alpha^{-1})= \min(x\lambda_{m,\hat{r}}^{-1})$ for all $x\in \{1\} \cup [p+2,m]$. Obviously, we have $[1,m-1]\beta = [1,m-1]\alpha = [1,m-1]\lambda_{m,\hat{r}} = 1$ and $(m-k)\beta^{-1} =(m-k)\alpha^{-1} =(m-k) \lambda_{m,\hat{r}}^{-1} =\{m+k\}$ for all $k \in [0,m-p-2]$. This shows that $\{1\}\cup [m,2m-p-2]$ is a transversal of both $\ker(\alpha)$ and $\ker(\beta)$. Let $\varepsilon$ be the unique idempotent in $\mathcal{OPD}(n,r)$ with $\im(\varepsilon) =\{1\}\cup [m,2m-p-2]$ and $\ker(\varepsilon)=\ker(\alpha)$. Obviously, $\beta= \varepsilon \alpha$. Therefore, we obtain $\beta\in \left\langle \mathcal{OPD}(n,r) ,\alpha \right\rangle$, and thus $L_{m,\hat{r}} \subseteq \left\langle \mathcal{OPD}(n,r),\alpha \right\rangle$. \qed

\begin{lemma}\label{l10}
	Let $3\leq m\leq n-\hat{r}+1$. Then $\mathcal{ORD}(n,r)\setminus L_{m,\hat{r}}$ is a maximal subsemigroup of $\mathcal{ORD}(n,r)$.
\end{lemma}

\proof Assume that there are $\beta, \gamma\in \mathcal{ORD}(n,r) \setminus L_{m,\hat{r}}$ such that $\beta \gamma\in L_{m,\hat{r}}$. Then $\im(\beta \gamma)= \im(\lambda_{m,\hat{r}}) =\{1\}\cup [p+2,m]$, where $p=\max\{0, m-\hat{r}\}$ and $\min(x(\beta \gamma)^{-1})= \min(x\lambda_{m,\hat{r}}^{-1})$ for all $x\in\{1\} \cup [p+2,m]$. Note that, $1,m \in \fix(\beta\gamma)$ implies $1,m \in \fix(\beta)\cap \fix(\gamma)$ since both $\beta$ and $\gamma$ are order-decreasing. Then, since $\beta\gamma|_{[m,2m-p-2]}$ is injective and monotone (order-reversing), it follows that $\beta_{\mid_{[m,2m-p-2]}}$ is also injective and monotone. 

If $\beta_{\mid_{[m,2m-p-2]}}$ is injective and order-preserving, then $m\beta=m$ implies that $\beta_{\mid_{[m,2m-p-2]}}$ is a partial identity and we obtain $\gamma_{\mid_{[m,2m-p-2]}} =\lambda_{m,\hat{r} \mid_{[m,2m-p-2]}}$. Since $1\gamma=1$, $m\gamma>(m+1)\gamma$, and $\gamma$ is oriented, we can conclude that $[1,m-1]\gamma \cap [2,m-1]=\emptyset$. Since $\gamma$ is order-decreasing, we get $[1,m-1]\gamma \subseteq [1,m-1]$, i.e. $[1,m-1]\gamma =\{1\}$. This shows $\gamma_{\mid_{[1,2m-p-2]}} =\lambda_{m,\hat{r} \mid_{[1,2m-p-2]}}$. Further, for $\gamma_{\mid_{ [2m-p-1,n]}}$, we have: If $p=0$, then $(2m-2)\gamma = 2$ and thus $[2m-1,n]\gamma =\{1\}$, i.e. $\gamma= \lambda_{m,\hat{r}}$ or $[2m-1,n]\gamma =\{2\}$ or there exists $s\in \{0,1,\ldots,n-2m\}$ such that $[2m-1,2m-1+s] \gamma =\{2\}$ and $[2m+s,n] \gamma =\{1\}$, since $\gamma$ is oriented. This implies $\gamma\in L_{m,\hat{r}}$, which is a contradiction. If $p=m-\hat{r}$, then $\hat{r}<m$ and $(m+\hat{r}-2) \gamma= m-\hat{r}+2 >2$. We observe that $\mathcal{RD}^{*}(n,r) \setminus \mathcal{RD}^{*}(n, \hat{r}) =\emptyset$ which implies that $\gamma\in \mathcal{RD}^{*}(n, \hat{r})$. Thus $[m+\hat{r}-1,n]\gamma = \{1\}$, i.e. $\gamma=\lambda_{m,\hat{r}}$ or $[m+\hat{r}-1,n] \gamma= \{m-\hat{r} +2\}$ or there exists $s \in \{0,1,\ldots,n-m-\hat{r}\}$ such that $[m+\hat{r}-1, m+\hat{r}-1+s] \gamma= \{m-\hat{r}+2\}$ and $[m+\hat{r}+s,n] \gamma =\{1\}$, since $\gamma$ is oriented. This implies $\gamma \in L_{m,\hat{r}}$, which is a contradiction.

Suppose now that $\beta_{\mid_{[m,2m-p-2]}}$ is order-reversing. Since $\gamma$ is order-decreasing and $\im(\beta \gamma_{\mid_{ [m,2m-p-2]}}) =\im(\lambda_{m, \hat{r}\mid_{ [m,2m-p-2]}}) =[p+2,m]$, we conclude that $\beta_{\mid_{ [m,2m-p-2]}}$ is an order-reversing bijection from $[m,2m-p-2]$ onto $[p+2,m]$ which shows that $\beta_{\mid_{ [m,2m-p-2]}} =\lambda_{m, \hat{r} \mid_{[m,2m-p-2]}}$, and by the same arguments as above, we obtain $\beta \in L_{m,\hat{r}}$, which is a contradiction.

Therefore, $\mathcal{ORD}(n,r)\setminus L_{m,\hat{r}}$ is a subsemigroup of $\mathcal{ORD}(n,r)$, which is maximal by Lemma \ref{l9}. \qed

\begin{theorem}\label{t8}
	Let $S$ be a subsemigroup of $\mathcal{ORD}(n,r)$. Then $S$ is maximal if and only if it belongs to one of the following types:
	\begin{itemize}
		\item $S=\mathcal{ORD}(n,r)\setminus\{\alpha\}$ for some $\alpha \in C$;
		\item $S=\mathcal{ORD}(n,r)\setminus\{\lambda_{m,\hat{r}}\}$ for some $n-\hat{r}+2\leq m\leq n-1$;
		\item $S=\mathcal{ORD}(n,r)\setminus L_{m,\hat{r}}$ for some $3\leq m\leq n-\hat{r}+1$.
	\end{itemize}
\end{theorem}

\proof Let $\alpha \in C$. By Proposition \ref{p4}, since $\alpha$ is undecomposable in $\mathcal{ORD}(n,r)$, we can conclude that $S= \mathcal{ORD}(n,r) \setminus \{\alpha\}$ is a maximal subsemigroup of $\mathcal{ORD}(n,r)$. Analogously, for $n-\hat{r}+2\leq m\leq n-1$ we have that $\lambda_{m, \hat{r}}$ is undecomposable in $\mathcal{ORD}(n,r)$ by Proposition \ref{p8}. Thus, we can conclude that $S= \mathcal{ORD}(n,r) \setminus \{\lambda_{m, \hat{r}}\}$ is a maximal subsemigroup of $\mathcal{ORD}(n,r)$ for $n-\hat{r}+2\leq m\leq n-1$. By Lemma \ref{l10}, we know that $S=\mathcal{ORD}(n,r) \setminus L_{m,\hat{r}}$ is a maximal subsemigroup of $\mathcal{ORD}(n,r)$ for $3\leq m<n-\hat{r}+2$.

Conversely, let $S$ be a maximal subsemigroup of $\mathcal{ORD}(n,r)$. Suppose that $C \nsubseteq S$. Then there is $\alpha \in C \setminus S$ and we obtain that $S=\mathcal{ORD}(n,r) \setminus \{\alpha\}$, since $\mathcal{ORD}(n,r) \setminus \{\alpha\}$ is a maximal subsemigroup of $\mathcal{ORD}(n,r)$ for all $\alpha \in C$.

Suppose that $C \subseteq S$. Since $C\cup G_{m,r}$ is a generating set of $\mathcal{ORD}(n,r)$ (see, Theorem \ref{t7}), we have $G_{m,r} \nsubseteq S$. Moreover, since $C \cup G_{m,r}$ is a minimal generating set of $\mathcal{ORD}(n,r)$ and $S$ is a maximal subsemigroup of $\mathcal{ORD}(n,r)$, we have $(\mathcal{ORD}(n,r)\setminus S) \cap G_{m,r} = \{\lambda_{m,\hat{r}}\}$ for some $3 \leq m \leq n-1$. If $n-\hat{r}+2\leq m\leq n-1$, then $S= \mathcal{ORD}(n,r) \setminus \{\lambda_{m, \hat{r}}\}$ since $\mathcal{ORD}(n,r) \setminus \{\lambda_{m, \hat{r}}\}$ is a maximal subsemigroup of $\mathcal{ORD}(n,r)$. Suppose that $3\leq m\leq n-\hat{r}+1$. Assume that $L_{m,\hat{r}} \cap S\neq \emptyset$. Then $L_{m,\hat{r}} \subseteq S$ by Lemma \ref{l9}, which contradicts with $\lambda_{m, \hat{r}}\notin S$. Therefore, $L_{m,\hat{r}} \cap S=\emptyset$, and so $S= \mathcal{ORD}(n,r) \setminus L_{m,\hat{r}}$, since $\mathcal{ORD}(n,r) \setminus L_{m,\hat{r}}$ is a maximal subsemigroup of $\mathcal{ORD}(n,r)$. \qed

Since $\mathcal{ORD}_{n} = \mathcal{ORD}(n,n-1) \cup \{1_{n}\}$, from Theorem \ref{t8}, we obtain the following result:

\begin{theorem}\label{t9}
	Let $S$ be a subsemigroup of $\mathcal{ORD}_n$. Then $S$ is a maximal subsemigroup of $\mathcal{ORD}_n$ if and only if $S = \mathcal{ORD}(n,n-1)$ or there exists a maximal subsemigroup $T$ of $\mathcal{ORD}(n,n-1)$ such that $S = T \cup \{1_n\}$.  \hfill $\square$
\end{theorem}

\end{document}